\documentclass[12pt,a4paper]{article}
\usepackage[T1]{fontenc}
\usepackage{lmodern}
\usepackage[margin=28mm]{geometry}
\usepackage{amsmath,amssymb,amsthm}
\usepackage{mathtools}
\usepackage{microtype}
\usepackage{tikz-cd}
\usepackage{cite}
\usepackage[hidelinks]{hyperref}
\hypersetup{
  pdftitle={Dominant root sequences of arbitrary degree},
  pdfauthor={Ling Chen, Chongrui Li and Minggang Wei},
  pdflang={en},pdfstartview={FitH}}

\theoremstyle{plain}
\newtheorem{theorem}{Theorem}[section]
\newtheorem{proposition}[theorem]{Proposition}
\newtheorem{lemma}[theorem]{Lemma}
\newtheorem{corollary}[theorem]{Corollary}
\theoremstyle{definition}
\newtheorem{example}[theorem]{Example}
\theoremstyle{remark}
\newtheorem{remark}[theorem]{Remark}
\numberwithin{equation}{section}

\newcommand{\C}{\mathbb C}
\newcommand{\A}{\mathbb A}
\newcommand{\GL}{\operatorname{GL}}
\newcommand{\SL}{\operatorname{SL}}
\newcommand{\End}{\operatorname{End}}
\newcommand{\diag}{\operatorname{diag}}
\newcommand{\id}{\operatorname{id}}

\title{\fontsize{18}{22}\selectfont\bfseries
Dominant root sequences of arbitrary degree}
\author{Ling Chen\qquad Chongrui Li\thanks{Corresponding author. E-mail:
\href{mailto:lichongrui24@mails.ucas.ac.cn}{\texttt{lichongrui24@mails.ucas.ac.cn}}.}\qquad Minggang Wei\\[0.75em]
{\normalsize School of Mathematical Sciences}\\
{\normalsize University of Chinese Academy of Sciences}\\
{\normalsize 19A Yuquan Road, Shijingshan District, Beijing 100049, China}\\[0.4em]
{\small\href{mailto:chenling2013@ucas.ac.cn}{\texttt{chenling2013@ucas.ac.cn}}\quad
\href{mailto:weiminggang21@mails.ucas.ac.cn}{\texttt{weiminggang21@mails.ucas.ac.cn}}}}
\date{}

\begin{document}
\maketitle
\begin{abstract}
\normalsize
We show that for root-subgroup multiplication maps, dominance does not
imply birationality. This gives a negative answer to a question of Fang,
Fourier and Littelmann. More precisely, for every integer $r\geq 2$, we
construct a sequence of positive roots in type $A_{3r-1}$; the associated
multiplication map is dominant of generic degree $r$, hence not
birational. In our construction,
a key observation is that Lusztig's flag incidence projection for a
triangle quiver is, up to birational equivalence and identity factors,
simultaneously a root-subgroup multiplication map and a multiplicative
partial Grothendieck--Springer map.
Moreover, we extend these examples to types $B$, $C$ and $D$
in every generic degree $r\geq 2$, and to types $E_6$, $E_7$ and $E_8$
in generic degree $2$.
\end{abstract}

\begingroup\small

\par\endgroup

\section{Introduction}

Toric degenerations of flag varieties translate representation-theoretic
and geometric invariants into polyhedral data: dimensions of
irreducible representations become lattice-point counts, and degrees
of projective embeddings become normalized volumes.
Two prominent sources of such degenerations are string parametrizations
of the dual canonical basis \cite{Caldero2002,AlexeevBrion2004} and the
PBW filtration together with its monomial refinements
\cite{FeiginFourierLittelmann2011,FeiginFourierLittelmann2017}.
The polytopes arising in both cases were realized as Newton--Okounkov
bodies: Kaveh obtained string polytopes from a highest-term valuation
associated with coordinates on a Bott--Samelson variety
\cite{Kaveh2015}; Kiritchenko obtained the
Feigin--Fourier--Littelmann--Vinberg polytopes in type~\(A\) from a
valuation defined by a flag of translated Schubert subvarieties
\cite{Kiritchenko2017}.

Fang, Fourier and Littelmann \cite{FFL2017} combined these viewpoints
in one framework. 
In their framework, the input is a sequence of positive roots, with length $N$
equal to the dimension of the open cell of the flag variety.
By multiplying the root subgroups of the positive roots in their order,
one can get a morphism from affine space to the open cell. 
This morphism is called the multiplication map of its corresponding root sequence.
When the multiplication map is birational, it gives rational coordinates on the open cell.
Moreover, with a suitable monomial order on \(\mathbb N^N\), these coordinates give a lowest-term valuation.
Furthermore, if the valuation monoid of this valuation is finitely
generated and saturated, one obtains normal toric degenerations of flag
and spherical varieties.
A sequence of positive roots is called birational if its multiplication map is birational \cite{FFL2017}.
Few birational
sequences were known beyond the examples in \cite{FFL2017, Bossinger2021}.

Fang, Fourier and Littelmann also gave a criterion for dominance in
\cite[Lemma~2]{FFL2017}. In the same paper, they raised a question
\cite[Remark~1]{FFL2017}, which can be restated as follows: does dominance of the
multiplication map imply birationality?
An affirmative answer would allow birational sequences to be recognized
by their dominance criterion alone.
Since the source and target have the same dimension, a dominant
multiplication map is generically finite. 
Thus their question is equivalent to whether every dominant
multiplication map has generic degree one.

In this paper, we show that the generic degree of a dominant multiplication map can be arbitrary, hence answering Fang, Fourier and Littelmann's question negatively.
Let \(r\) be any positive integer. For any sequence \(r_1,\dots,r_s\) of positive integers with \(\sum_{i=1}^{s}r_i= r\), we construct a dominant multiplication map of generic degree \(r!/(r_1!\cdots r_s!)\) in type \(A_{3r-1}\).
Taking \(r_1=1\) and \(\ r_2=r-1\) with \(r\ge 2\), one immediately sees that dominant multiplication maps can be of arbitrary degree.
In the case \(r_1=r_2=\cdots=r_s=1\), we also prove that the associated function-field extension is Galois with symmetric group \(\mathfrak S_r\). 
Moreover, by passing to Levi subgroups, we obtain dominant multiplication maps of degree \(r\ (r\ge 2)\) in types
\(B\), \(C\) and \(D\), and of degree \(2\) in types \(E_6\), \(E_7\) and \(E_8\).
In these constructions, our key observation is that Lusztig's flag incidence projection \cite{Lusztig1991}
for a triangle quiver is, up to birational equivalence and identity factors, both a root-subgroup
multiplication map and a multiplicative partial Grothendieck--Springer map.
Since birational equivalence preserves generic degree, the second description
determines the degree of the first.

Throughout this paper, let \(G\) be a connected semisimple algebraic group over the complex field
\(\mathbb C\).
Fix a maximal torus \(T\subset G\) and a system \(\Phi^+\) of positive
roots in the root system \(\Phi\) of \((G,T)\), and let \(U^-\) be the
corresponding negative maximal unipotent subgroup, which is identified
with the open cell of the flag variety of \(G\); set
\(N=|\Phi^+|=\dim U^-\).
For each \(\beta\in\Phi^+\), fix an isomorphism
\(x_{-\beta}:\mathbb G_a\xrightarrow{\sim}U_{-\beta}\) from the
additive group onto the negative root subgroup
\(U_{-\beta}\subset U^-\).
For a sequence \(S=(\beta_1,\ldots,\beta_N)\) of positive roots,
with repetitions allowed, the \emph{root-subgroup multiplication map} is
\begin{equation*}
\pi_S:\mathbb A^N\longrightarrow U^-,
\qquad
(t_1,\ldots,t_N)\longmapsto
x_{-\beta_1}(t_1)\cdots x_{-\beta_N}(t_N).
\end{equation*}
Following \cite[Definition~2]{FFL2017}, we call \(S\)
\emph{birational} if \(\pi_S\) is birational.
Moreover, for convenience, we shall also call \(S\) \emph{dominant} if \(\pi_S\)
is dominant.
For a dominant sequence \(S\), we write
\(\deg\pi_S=[\mathbb C(\mathbb A^N):\pi_S^*\mathbb C(U^-)]\) for the
generic degree of \(\pi_S\).
Any other choice of isomorphisms \(x_{-\beta}\) does not affect these notions, since it only rescales each coordinate
\(t_i\) by a nonzero constant.

Let \(f:X\to Y\) and \(f':X'\to Y'\) be dominant morphisms of irreducible varieties.
We say that \(f\) and \(f'\) are \emph{birationally equivalent}, and write
\(f\sim_{\mathrm{bir}}f'\), if there are birational maps
\(\alpha:X\dashrightarrow X'\) and \(\beta:Y\dashrightarrow Y'\) satisfying
\(\beta\circ f=f'\circ\alpha\).
Such maps \(\alpha\) and \(\beta\) identify the function-field extensions of \(f\) and \(f'\).
Hence birationally equivalent morphisms have the same generic degree, and also the same Galois group if the extensions are Galois.

Next we recall 
the necessary notions about the
multiplicative partial Grothendieck--Springer maps.
Let \(\boldsymbol{\vec{r}}=(r_1,\dots,r_s)\) be a composition of \(r\), that is, a sequence of positive integers with \(\sum_{i=1}^{s}r_i= r\).
Let \(P_{\boldsymbol{\vec{r}}}\subset\GL_r\) be the standard parabolic subgroup
with diagonal blocks of sizes \(r_1,\ldots,r_s\), and let
\(\mathcal F_{\boldsymbol{\vec{r}}}=\GL_r/P_{\boldsymbol{\vec{r}}}\) be the variety of flags
\(F_\bullet:0=F_0\subset F_1\subset\cdots\subset F_s=\mathbb C^r\)
with \(\dim F_i=r_1+\cdots+r_i\).
Set
\begin{equation}\label{eq:flag-cover}
X_{\boldsymbol{\vec{r}}}
=\{(A,F_\bullet)\in\GL_r\times\mathcal F_{\boldsymbol{\vec{r}}}:
AF_i=F_i\text{ for all }i\}.
\end{equation}
The \emph{multiplicative partial Grothendieck--Springer map} is the
projection 
\begin{equation}\label{eq:flag-proj}
  p_{\boldsymbol{\vec{r}}}:X_{\boldsymbol{\vec{r}}}\to\GL_r,\qquad
(A,F_\bullet)\mapsto A.
\end{equation}

Then we present the main result of this paper in the following theorem:
\begin{theorem}\label{thm:flags}
Let \(r\ge1\) and let \(\boldsymbol{\vec{r}}=(r_1,\ldots,r_s)\) be a composition
of \(r\).
Then there exists a dominant root sequence \(S_{\boldsymbol{\vec{r}}}\) of length
\(\binom{3r}{2}\) in type \(A_{3r-1}\), such that
\begin{equation}\label{eq:flag-normal-form}
\pi_{S_{\boldsymbol{\vec{r}}}}
\sim_{\mathrm{bir}}
\id_{\A^{N-r^2}}\times p_{\boldsymbol{\vec{r}}}.
\end{equation}
Consequently, \(\deg\pi_{S_{\boldsymbol{\vec{r}}}}=r!/(r_1!\cdots r_s!)\), and
for \(\boldsymbol{\vec{r}}=(1,\ldots,1)\), the extension \(\C(\A^N)/\pi_{S_{\boldsymbol{\vec{r}}}}^*\C(U^-)\)
is Galois with
symmetric group \(\mathfrak S_r\).
\end{theorem}

Taking \(\boldsymbol{\vec r}=(1,r-1)\) with \(r\ge2\) in the above theorem, we get \(\deg\pi_{S_{\boldsymbol{\vec r}}}=r\ge 2\).
Thus we obtain dominant root sequences of arbitrary generic degree.
Moreover, the proof for Theorem~\ref{thm:flags} in Section~\ref{sec:three-blocks} gives an explicit construction of the dominant root sequence \(S_{\boldsymbol{\vec{r}}}\).
We show such a construction for $A_5$ in Example~\ref{ex:A5}.
The resulting root sequence, which is of generic degree $2$, is as follows:
\begin{equation}\label{eq:A5-example}
S=\bigl(
\alpha_1,\alpha_3,\alpha_5,
[1,2],[1,4],
\alpha_2,[2,3],[2,4],[2,5],
[3,4],\alpha_4,[4,5],
\alpha_1,\alpha_3,\alpha_5
\bigr),
\end{equation}
where
\(\alpha_1,\ldots,\alpha_5\) are the simple roots of \(A_5\) in the standard order, 
and
\([i,j]=\alpha_i+\cdots+\alpha_j\) for \(1\le i<j\le 5\).

This paper is organized as follows.
In Section~\ref{sec:flag-cover}, we recall some geometric properties of the multiplicative partial Grothendieck--Springer map 
\(p_{\boldsymbol{\vec{r}}}\), and relate \(p_{\boldsymbol{\vec{r}}}\) to Lusztig's incidence
projections for some quivers.
In Section~\ref{sec:three-blocks}, we realize incidence projections of
acyclic quivers by root-subgroup multiplications, construct  dominant root sequences
\(S_{\boldsymbol{\vec{r}}}\) for type \(A\), and hence prove Theorem~\ref{thm:flags}.
In Section~\ref{sec:levi}, we extend the sequences \(S_{\boldsymbol{\vec{r}}}\) for type \(A\) to types \(B\), \(C\), \(D\) and \(E\).

\section{Invariant flags and quiver incidence}
\label{sec:flag-cover}

In this section, we first
recall some geometric  properties of 
the multiplicative partial Grothendieck--Springer map \(p_{\boldsymbol{\vec{r}}}\) defined in 
\eqref{eq:flag-proj}.
Then we introduce quiver flag incidence projections and relate a class
of these projections to \(p_{\boldsymbol{\vec{r}}}\).

We begin with the geometry of \(X_{\boldsymbol{\vec{r}}}\).
Let \(F_\bullet^{\boldsymbol{\vec{r}}}\) be the standard partial flag of type
\(\boldsymbol{\vec{r}}=(r_1,\ldots,r_s)\), so that
\(P_{\boldsymbol{\vec{r}}}=\operatorname{Stab}_{\GL_r}(F_\bullet^{\boldsymbol{\vec{r}}})\).
The associated-bundle description gives a canonical isomorphism
\begin{equation}\label{eq:associated-bundle}
\GL_r\times^{P_{\boldsymbol{\vec{r}}}}P_{\boldsymbol{\vec{r}}}
\xrightarrow{\sim}X_{\boldsymbol{\vec{r}}},
\qquad
[g,h]\longmapsto(ghg^{-1},gF_\bullet^{\boldsymbol{\vec{r}}}),
\end{equation}
where \(x\in P_{\boldsymbol{\vec{r}}}\) acts by
\((g,h)\mapsto(gx,x^{-1}hx)\).
Hence, \(X_{\boldsymbol{\vec{r}}}\) is irreducible and
\(\dim X_{\boldsymbol{\vec{r}}}=\dim\mathcal F_{\boldsymbol{\vec{r}}}+\dim P_{\boldsymbol{\vec{r}}}=r^2\).

Let \(\GL_r^{\mathrm{rs}}\subset\GL_r\) denote the regular semisimple
locus, consisting of matrices with \(r\) distinct eigenvalues. Then \(\GL_r^{\mathrm{rs}}\) is a
dense open subset of \(\GL_r\).
For convenience of the reader, we include a short proof of the standard facts
on the generic degree of \(p_{\boldsymbol{\vec{r}}}\) in (\ref{eq:flag-proj})
and the corresponding extension of function fields in the following.

\begin{proposition}\label{prop:flag-cover}
The multiplicative partial Grothendieck--Springer map \(p_{\boldsymbol{\vec{r}}}\) defined in (\ref{eq:flag-proj}) is dominant, with generic degree
\({r!}/({r_1!\cdots r_s!})\).
For \(\boldsymbol{\vec{r}}=(1,\ldots,1)\),
the extension \(\C(X_{\boldsymbol{\vec{r}}})/p_{\boldsymbol{\vec{r}}}^*\C(\GL_r)\) is 
Galois with symmetric group \(\mathfrak S_r\).
\end{proposition}

\begin{proof} Since \(\dim X_{\boldsymbol{\vec{r}}}=\dim\GL_r\), the generic degree of \(p_{\boldsymbol{\vec{r}}}\)
is the number of points in a general fiber.
Hence to compute the generic degree of \(p_{\boldsymbol{\vec{r}}}\), it suffices to compute \(|p_{\boldsymbol{\vec{r}}}^{-1}(A)|\) for every \(A\) in the dense
open subset \(\GL_r^{\mathrm{rs}}\).
For $A\in\GL_r^{\mathrm{rs}}$, write
$\C^r=L_1\oplus\cdots\oplus L_r$ for its decomposition into eigenlines.
Since the eigenvalues are distinct, every $A$-invariant subspace is a
direct sum of some of the $L_j$'s.
Hence the $A$-invariant flags in $\mathcal F_{\boldsymbol{\vec{r}}}$ correspond to
the ordered partitions $(I_1,\ldots,I_s)$ of $\{1,\ldots,r\}$ with
$|I_i|=r_i$, and there are $r!/(r_1!\cdots r_s!)$ such partitions.

For $\boldsymbol{\vec{r}}=(1,\ldots,1)$, the invariant complete flags of
$A\in\GL_r^{\mathrm{rs}}$ correspond to the orderings of the eigenlines.
Reordering them gives an algebraic action of $\mathfrak S_r$ on
$p_{\boldsymbol{\vec{r}}}^{-1}(\GL_r^{\mathrm{rs}})$ that commutes with
$p_{\boldsymbol{\vec{r}}}$ and is simply transitive on each fiber.
Since $p_{\boldsymbol{\vec{r}}}^{-1}(\GL_r^{\mathrm{rs}})$ is dense in
$X_{\boldsymbol{\vec{r}}}$, this gives a faithful action of $\mathfrak S_r$ on
$\C(X_{\boldsymbol{\vec{r}}})$ by $p_{\boldsymbol{\vec{r}}}^*\C(\GL_r)$-automorphisms.
As $[\C(X_{\boldsymbol{\vec{r}}}):p_{\boldsymbol{\vec{r}}}^*\C(\GL_r)]=r!=|\mathfrak S_r|$,
the extension is Galois with group $\mathfrak S_r$.
\end{proof}

To connect \(p_{\boldsymbol{\vec{r}}}\) with root-subgroup multiplication, we now
consider quiver representations, following Lusztig
\cite[Sections~1.3--1.6]{Lusztig1991}.
We first introduce the representation space \(R_Q\) and the graded flag
variety \(\mathcal F_{\boldsymbol\lambda}\).

Let \(Q=(Q_0,Q_1)\) be a finite quiver with vertex set \(Q_0\) and
arrow set \(Q_1\), and write \(|Q|\) for its underlying undirected
graph. Throughout the rest of the paper, we assume that \(Q\) is simple, that is, \(|Q|\) is a simple graph.
Fix a finite-dimensional \(Q_0\)-graded complex vector space
\(W=\bigoplus_{i\in Q_0}W_i\). The representations of \(Q\) on
\(W\) form the affine space
\begin{equation}\label{eq:quiver-representation-space}
R_Q=\bigoplus_{(\gamma:j\to i)\in Q_1}\operatorname{Hom}(W_j,W_i),
\end{equation}
whose points are written as \(A=(A_\gamma)_{\gamma\in Q_1}\).

For a chain \(F_\bullet:0=F_0\subseteq F_1\subseteq\cdots\subseteq
F_s=W\) of subspaces, set \(F_a^{(i)}=F_a\cap W_i\) and
\(F_\bullet^{(i)}=(F_a^{(i)})_{a=0}^s\).
The chain \(F_\bullet\) is a \emph{\(Q_0\)-graded flag} if
\(F_a=\bigoplus_{i\in Q_0}F_a^{(i)}\) for all \(a\).
Its \emph{type} is the sequence
\(\boldsymbol\lambda=(\lambda^1,\ldots,\lambda^s)\) in
\(\mathbb Z_{\ge0}^{Q_0}\) given by
\(\lambda^a=(\dim (F_a^{(i)}/F_{a-1}^{(i)}))_{i\in Q_0}\).
Fix a sequence \(\boldsymbol\lambda=(\lambda^1,\ldots,\lambda^s)\) in
\(\mathbb Z_{\ge0}^{Q_0}\) with
\(\lambda^1+\cdots+\lambda^s=(\dim W_i)_{i\in Q_0}\), and write
\(\lambda^a=(\lambda_i^a)_{i\in Q_0}\).
Let \(\mathcal F_{\boldsymbol\lambda}\) be the variety of
\(Q_0\)-graded flags in \(W\) of type \(\boldsymbol\lambda\).

For \(A\in R_Q\), a flag \(F_\bullet\in\mathcal F_{\boldsymbol\lambda}\) is called \emph{\(A\)-stable} if
\(A_\gamma F_a^{(j)}\subseteq F_a^{(i)}\)
for every \(1\le a<s\) and every arrow \(\gamma:j\to i\).
Set
\begin{equation}\label{eq:quiver-incidence}
I_{Q,\boldsymbol\lambda}
=\bigl\{(A,F_\bullet)\in
R_Q\times\mathcal F_{\boldsymbol\lambda}:
F_\bullet\text{ is }A\text{-stable}\bigr\}.
\end{equation}
The \emph{incidence projection} is the map
\begin{equation}\label{eq:quiver-proj}
\rho_{Q,\boldsymbol\lambda}:I_{Q,\boldsymbol\lambda}\to R_Q,\qquad
(A,F_\bullet)\mapsto A.
\end{equation}

We next find a necessary condition for \(\rho_{Q,\boldsymbol\lambda}\)
to be dominant and generically finite.
Such a morphism requires \(\dim I_{Q,\boldsymbol\lambda}=\dim R_Q\), so we first
compute \(\dim I_{Q,\boldsymbol\lambda}-\dim R_Q\) for arbitrary \(Q\) and \(\boldsymbol{\lambda}\).

Let the group \(G_W=\prod_{i\in Q_0}\GL(W_i)\) act on
\(\mathcal F_{\boldsymbol\lambda}\) naturally, and on \(R_Q\) by
\((g.A)_\gamma=g_iA_\gamma g_j^{-1}\) for \(\gamma:j\to i\), where \(g\in G_W\), \(A\in R_Q\).
Fix an ordered basis for each \(W_i\), and let
\(F_\bullet^{\boldsymbol\lambda}\in\mathcal F_{\boldsymbol\lambda}\) be
the graded flag such that \(F_a^{\boldsymbol\lambda}\cap W_i\) is
spanned by the first \(\lambda_i^1+\cdots+\lambda_i^a\) vectors  of the fixed basis of
\(W_i\).
Since \(G_W\) acts transitively on \(\mathcal F_{\boldsymbol\lambda}\), the stabilizer
\(P=\operatorname{Stab}_{G_W}(F_\bullet^{\boldsymbol\lambda})\) gives
\(\mathcal F_{\boldsymbol\lambda}\simeq G_W/P\).
Let \(Z=\{B\in R_Q:(B,F_\bullet^{\boldsymbol\lambda})\in
I_{Q,\boldsymbol\lambda}\}\); \(Z\) is a linear subspace of \(R_Q\).
For \(\mathbf u,\mathbf v\in\mathbb Z^{Q_0}\), let
\(\langle\mathbf u,\mathbf v\rangle_Q
=\sum_{i\in Q_0}u_iv_i-\sum_{\gamma:j\to i}u_jv_i\)
be the \emph{Euler form} of \(Q\).
The following lemma expresses
\(\dim I_{Q,\boldsymbol\lambda}-\dim R_Q\) in terms of the Euler form.

\begin{lemma}\label{lem:quiver-incidence-dimension}
The variety \(I_{Q,\boldsymbol\lambda}\)   defined in (\ref{eq:quiver-incidence}) is irreducible, and
\begin{equation}\label{eq:quiver-dimension}
\dim I_{Q,\boldsymbol\lambda}-\dim R_Q
=\sum_{a<b}\langle\lambda^a,\lambda^b\rangle_Q.
\end{equation}
\end{lemma}

\begin{proof}Since \(P\) fixes \(F_\bullet^{\boldsymbol\lambda}\), the subspace \(Z\) is \(P\)-stable.
The projection \(I_{Q,\boldsymbol\lambda}\to\mathcal F_{\boldsymbol\lambda}\simeq G_W/P\)
is \(G_W\)-equivariant with fiber \(Z\) over \(F_\bullet^{\boldsymbol\lambda}\), so
there is an isomorphism
\begin{equation}\label{eq:quiver-bundle}
G_W\times^P Z\xrightarrow{\sim}I_{Q,\boldsymbol\lambda},
\qquad
[g,B]\longmapsto\bigl(g.B,\,gF_\bullet^{\boldsymbol\lambda}\bigr),
\end{equation}
where \(x\in P\) acts on \(G_W\times Z\) by \((g,B)\mapsto(gx,x^{-1}.B)\).
Thus \(I_{Q,\boldsymbol\lambda}\) is a vector bundle over \(G_W/P\), hence irreducible.

By \eqref{eq:quiver-bundle},
\(\dim I_{Q,\boldsymbol\lambda}=\dim(G_W/P)+\dim Z\), so
\[
\dim I_{Q,\boldsymbol\lambda}-\dim R_Q
=\dim(G_W/P)-\operatorname{codim}_{R_Q}Z.
\]
Since \(G_W/P\) is a product of partial flag varieties,
\(\dim(G_W/P)=\sum_{i\in Q_0}\sum_{a<b}\lambda_i^a\lambda_i^b\).

Write \(F_k^{\boldsymbol\lambda,(i)}=F_k^{\boldsymbol\lambda}\cap W_i\).
For \(i\in Q_0\) and \(1\le a\le s\), let \(W_i^a\) be the span of the basis vectors of
\(W_i\) lying in \(F_a^{\boldsymbol\lambda,(i)}\) but not in \(F_{a-1}^{\boldsymbol\lambda,(i)}\).
Then \(\dim W_i^a=\lambda_i^a\) and \(F_k^{\boldsymbol\lambda,(i)}=\bigoplus_{a\le k}W_i^a\).
For \(\gamma:j\to i\), let \(B_\gamma^{ba}:W_j^a\to W_i^b\) be the \((b,a)\)-block of \(B_\gamma\).
Then \(B_\gamma F_k^{\boldsymbol\lambda,(j)}\subseteq F_k^{\boldsymbol\lambda,(i)}\) for all \(k\)
if and only if \(B_\gamma^{ba}=0\) for all \(b>a\).
Hence
\(\operatorname{codim}_{R_Q}Z=\sum_{\gamma:j\to i}\sum_{a<b}\lambda_j^a\lambda_i^b\).
\end{proof}

Now we relate a class of \(\rho_{Q,\boldsymbol\lambda}\) to \(p_{\boldsymbol{\vec{r}}}\).
Fix a composition \(\boldsymbol{\vec{r}}=(r_1,\ldots,r_s)\) of \(r\),
and take \(\dim W_i=r\) and \(\lambda_i^a=r_a\) for \(i\in Q_0\) and \(1\le a\le s\). 
Noting that \(\rho_{Q,\boldsymbol\lambda}\) is the product of the incidence
projections of the connected components of \(Q\), for simplicity, we may assume
\(Q\) connected.
Let \(\#(X)\) denote the number of elements in the set \(X\).
Since \(\langle\lambda^a,\lambda^b\rangle_Q=r_ar_b(\# (Q_0)-\#(Q_1))\) and
\(\sum_{a<b}r_ar_b=\dim\mathcal F_{\boldsymbol{\vec{r}}}\),
Eq.~\eqref{eq:quiver-dimension} becomes
\(\dim I_{Q,\boldsymbol\lambda}-\dim R_Q
=(\#(Q_0)-\#(Q_1))\dim\mathcal F_{\boldsymbol{\vec{r}}}\). 
Hence, for  \(s\ge2\), \(\dim I_{Q,\boldsymbol\lambda}=\dim R_Q\) if and only if
\(\# (Q_0)=\#(Q_1)\).
Moreover, \(\# (Q_0)=\#(Q_1)\) is equivalent to \(|Q|\) being unicyclic,
i.e., \(|Q|\) has exactly one cycle. The following lemma shows that, when \(|Q|\) is unicyclic, 
\(\rho_{Q,\boldsymbol\lambda}\) can be reduced to \(p_{\boldsymbol{\vec{r}}}\).

\begin{lemma}\label{lem:unicyclic-flag-projection}
Let \(Q=(Q_0,Q_1)\) be a connected simple quiver such that \(|Q|\) is unicyclic.
Let \(\boldsymbol{\vec{r}}=(r_1,\ldots,r_s)\) be a composition of \(r\). When
\(\dim W_i=r\) and \(\lambda_i^a=r_a\) for all \(i\in Q_0\) and \(1\le a\le s\),
the incidence projection \(\rho_{Q,\boldsymbol\lambda}\) defined in (\ref{eq:quiver-proj}) is dominant, and
\begin{equation*}
\rho_{Q,\boldsymbol\lambda}\sim_{\mathrm{bir}}
\id_{(\GL_r)^{\#(Q_0)-1}}\times p_{\boldsymbol{\vec{r}}}.
\end{equation*}
\end{lemma}

\begin{proof}

Fix identifications \(W_i\simeq\C^r\) for \(i\in Q_0\), and set
\(R^\times=\{A\in R_Q:\det A_\gamma\ne0\text{ for all }\gamma\in Q_1\}\)
and \(I^\times=\rho_{Q,\boldsymbol\lambda}^{-1}(R^\times)\).
Then \(R^\times\) is dense open in \(R_Q\), and \(I^\times\) is open in
\(I_{Q,\boldsymbol\lambda}\).
Observe that \((A,F_\bullet)\) lies in \(I^\times\) when \(A_\gamma=I_r\), \(F_\bullet^{(i)}=F_\bullet^{\boldsymbol{\vec{r}}}\) for all \(\gamma\in Q_1\), \(i\in Q_0\).
Hence \(I^\times\) is nonempty.
Since \(I_{Q,\boldsymbol\lambda}\) is irreducible by
Lemma~\ref{lem:quiver-incidence-dimension}, \(I^\times\) is dense.

Choose a spanning tree \(\mathcal T\) of \(|Q|\); since \(|Q|\) is unicyclic,
exactly one arrow \(\gamma_0:u\to v\) of \(Q\) lies outside \(\mathcal T\).
For \(A\in R^\times\) and \(i\in Q_0\), let \(v=i_0,i_1,\ldots,i_m=i\) be the path
in \(\mathcal T\) from \(v\) to \(i\), and let \(\gamma_k\) be the arrow of \(\mathcal T\)
between \(i_{k-1}\) and \(i_k\). Set
\[
g_i=A_{\gamma_m}^{\varepsilon_m}\cdots A_{\gamma_1}^{\varepsilon_1}\in\GL_r,
\qquad
\varepsilon_k=
\begin{cases}
1 & \text{if } \gamma_k:i_{k-1}\to i_k,\\
-1 & \text{if } \gamma_k:i_k\to i_{k-1}.
\end{cases}
\]
Then \(g_v=I_r\) and \(g_i=A_\gamma g_j\) for every arrow \(\gamma:j\to i\) in \(\mathcal T\).
Since all \(A_\gamma\)'s are invertible and all \(F_\bullet^{(k)}\)'s have type
\(\boldsymbol{\vec{r}}\), \(F_\bullet\) is stable along the arrow \(\gamma:j\to i\)
if and only if \(A_\gamma F_\bullet^{(j)}=F_\bullet^{(i)}\).
Hence stability along \(\mathcal T\) is equivalent to
\(F_\bullet^{(i)}=g_iF_\bullet^{(v)}\) for all \(i\in Q_0\). Therefore, stability
along \(\gamma_0\) is equivalent to
\((A_{\gamma_0}g_u,F_\bullet^{(v)})\in X_{\boldsymbol{\vec{r}}}\).

After ordering the tree arrows, we only need to prove that the following diagram commutes and that its horizontal maps are isomorphisms:
\begin{equation*}
\begin{tikzcd}[column sep=large]
I^\times \arrow[r,"\phi","\sim"']
\arrow[d,"\rho_{Q,\boldsymbol\lambda}"']
&(\GL_r)^{\# (Q_0)-1}\times X_{\boldsymbol{\vec{r}}}
\arrow[d,"\id\times p_{\boldsymbol{\vec{r}}}"]\\
R^\times \arrow[r,"\psi","\sim"']
&(\GL_r)^{\# (Q_0)-1}\times\GL_r,
\end{tikzcd}
\end{equation*}
where
\begin{equation}\label{eq:free-factor}
\phi(A,F_\bullet)
=\bigl((A_\gamma)_{\gamma\in\mathcal T},
(A_{\gamma_0}g_u,F_\bullet^{(v)})\bigr),\qquad
\psi(A)
=\bigl((A_\gamma)_{\gamma\in\mathcal T},
A_{\gamma_0}g_u\bigr).
\end{equation}
Note that each \(g_i\) depends only on \((A_\gamma)_{\gamma\in\mathcal T}\).
Let \(\bigl((B_\gamma)_{\gamma\in\mathcal T},(C,F_\bullet)\bigr)\in
(\GL_r)^{\#(Q_0)-1}\times X_{\boldsymbol{\vec{r}}}\), define \(g_i\) as above using the
\(B_\gamma\), and let \((A,F'_\bullet)\) be given by
\[
A_\gamma=B_\gamma\ (\gamma\in\mathcal T),\qquad A_{\gamma_0}=Cg_u^{-1},\qquad
F'^{(i)}_\bullet=g_iF_\bullet\ (i\in Q_0).
\]
By the previous paragraph, \((A,F'_\bullet)\in I^\times\) and
\(\phi(A,F'_\bullet)=\bigl((B_\gamma)_{\gamma\in\mathcal T},(C,F_\bullet)\bigr)\).
Conversely, every \((A,F_\bullet)\in I^\times\) arises in this way from
\(\phi(A,F_\bullet)\), since \(F_\bullet^{(i)}=g_iF_\bullet^{(v)}\).
Hence \(\phi\) is bijective, and similarly \(\psi\) is bijective with inverse
\(\bigl((B_\gamma)_{\gamma\in\mathcal T},C\bigr)\mapsto A\).
Since each \(g_i\) is a product of \(A_\gamma^{\pm1}\)'s, these maps are regular, and
\(\psi\circ\rho_{Q,\boldsymbol\lambda}
=(\id\times p_{\boldsymbol{\vec{r}}})\circ\phi\).
Since \(p_{\boldsymbol{\vec{r}}}\) is dominant by
Proposition~\ref{prop:flag-cover}, the diagram proves both
assertions.
\end{proof}

\section{Root-subgroup realization}
\label{sec:three-blocks}

In this section, we show that for every acyclic quiver \(Q\) and flag type \(\boldsymbol\lambda\), the restriction of
\(\rho_{Q,\boldsymbol\lambda}\) to a dense open subset of
\(I_{Q,\boldsymbol\lambda}\) is, up to isomorphism and an identity
factor, a root-subgroup multiplication map.
We first write \(\rho_{Q,\boldsymbol\lambda}\) as the multiplication
map of three affine subvarieties of a lower unitriangular group, and
then realize the three affine subvarieties by root subgroups.
With these results, we subsequently give the proof of Theorem~\ref{thm:flags}.

To place \(R_Q\) inside a lower unitriangular group, we number the
vertices in \(Q_0\) by \(\{1,\ldots,n\}\) so that \(j<i\) for every arrow
\(j\to i\) in \(Q_1\); such a numbering exists if and only if \(Q\) is acyclic.
In the rest of this section, we always assume \(Q\) is acyclic with such a numbering.
We order the basis of \(W\) by concatenating the fixed bases of
\(W_1,\ldots,W_n\). Take \(G=\SL(W)\), with \(T\)
the diagonal torus and \(\Phi^+\) the roots of the upper triangular
Borel subgroup.
Thus \(U^-\subset\SL(W)\) is the lower unitriangular subgroup, and
\(N=\binom{\dim W}{2}\).
Since \(Q\) is simple, we regard \(R_Q\) as a subspace of \(\End(W)\) by
identifying \(A\in R_Q\) with the block matrix \((A_{ij})_{i,j\in Q_0}\), where \(A_{ij}\) is the matrix of \(A_\gamma\) in the fixed bases if
\(\gamma:j\to i\), and \(A_{ij}=0\) otherwise.
As \(j<i\) for every arrow \(j\to i\), we have \(I_{\dim W}+A\in U^-\), where \(I_{\dim W}\) is the identity matrix and \(A\in R_Q\).

To relate \(U^-\) to the target \(R_Q\) of \(\rho_{Q,\boldsymbol\lambda}\),
we write elements of \(U^-\) as block matrices with respect to
\(W=\bigoplus_{i\in Q_0}W_i\).
Recall that \(G_W=\prod_{i\in Q_0}\GL(W_i)\). Set \(H=G_W\cap U^-\).
For \(g\in U^-\), write \(g_{ij}\in\operatorname{Hom}(W_j,W_i)\) for the
\((i,j)\)-block of \(g\).
Since \(Q\) is simple and \(j<i\) for every arrow \(j\to i\), the vector space of strictly lower blocks splits as
\begin{equation*}
\bigoplus_{i>j}\operatorname{Hom}(W_j,W_i)=R_Q\oplus R_Q',
\qquad\text{where}\quad
R_Q'=\bigoplus_{\substack{i>j\\ \text{no arrow }j\to i}}
\operatorname{Hom}(W_j,W_i).
\end{equation*}
For \(g\in U^-\), let \(h=\diag(g_{11},\ldots,g_{nn})\in H\).
Then \(gh^{-1}\) has identity diagonal blocks. 
Thus
\(gh^{-1}=I_{\dim W}+B+B'\) with \(B\in R_Q\) and \(B'\in R_Q'\).
Conversely, if \(g=(I_{\dim W}+B+B')h'\) with \(B\in R_Q\), \(B'\in R_Q'\) and
\(h'\in H\), then comparing diagonal blocks gives \(h'=h\).
Hence we get an isomorphism
\begin{equation}\label{eq:mu}
\mu:R_Q\times R_Q'\times H\xrightarrow{\sim}U^-,
\qquad (B,B',h)\longmapsto(I_{\dim W}+B+B')h.
\end{equation}
Define \(q:U^-\to R_Q\) by \(q(g)_\gamma=g_{ij}g_{jj}^{-1}\) for each
arrow \(\gamma:j\to i\).
Then \(q\) is the \(R_Q\)-coordinate of \(\mu^{-1}\).
Since \(R_Q'\) and \(H\) are affine spaces, \(\mu\) gives an
isomorphism over \(R_Q\):
\begin{equation}\label{eq:quiver-target-product}
U^-\simeq\A^{N-\dim R_Q}\times R_Q.
\end{equation}

Next, we consider the matrices in \(U^-\) with identity diagonal blocks
and arrow blocks keeping \(F_\bullet^{\boldsymbol\lambda}\) stable. Let
\begin{equation}\label{eq:quiver-standard-matrices}
Y_{Q,\boldsymbol\lambda}
=\{y\in U^-:y_{ii}=I_{\dim W_i}\text{ for all }i \in Q_0,\
(q(y),F_\bullet^{\boldsymbol\lambda})\in I_{Q,\boldsymbol\lambda}\}.
\end{equation}
Then \(\mu^{-1}(Y_{Q,\boldsymbol\lambda})=Z\times R_Q'\times\{I_{\dim W}\}\), hence 
\(Y_{Q,\boldsymbol\lambda}\) is an affine space.

Now we restrict \(\rho_{Q,\boldsymbol\lambda}\) to flags in an open
cell of \(\mathcal F_{\boldsymbol\lambda}\).
Let \(V\subset H\) be the unipotent radical of the parabolic subgroup
of \(G_W\) opposite to \(P\).
With respect to the above fixed basis of \(W\), we have \(V=\prod_{i\in Q_0}V_i\), where \(V_i\subset\GL(W_i)\)
consists of the block lower unitriangular matrices with diagonal blocks of
sizes \(\lambda_i^1,\ldots,\lambda_i^s\).
Then \(\mathcal F_{\boldsymbol\lambda}^{\circ}
=\{vF_\bullet^{\boldsymbol\lambda}:v\in V\}\) is open in
\(\mathcal F_{\boldsymbol\lambda}\), and the orbit map in the following is isomorphic:
\begin{equation}\label{eq:quiver-chart}
V\xrightarrow{\sim}\mathcal F_{\boldsymbol\lambda}^{\circ},
\qquad v\longmapsto vF_\bullet^{\boldsymbol\lambda}.
\end{equation}
Set \(I_{Q,\boldsymbol\lambda}^{\circ}=I_{Q,\boldsymbol\lambda}\cap
(R_Q\times\mathcal F_{\boldsymbol\lambda}^{\circ})\).
Since \((0,F_\bullet^{\boldsymbol\lambda})\in
I_{Q,\boldsymbol\lambda}^{\circ}\) and \(I_{Q,\boldsymbol\lambda}\) is
irreducible by Lemma~\ref{lem:quiver-incidence-dimension},
\(I_{Q,\boldsymbol\lambda}^{\circ}\) is dense open in
\(I_{Q,\boldsymbol\lambda}\).

The morphisms \(f:X\to Y\) and \(f':X'\to Y'\) are called
\emph{isomorphic}, and write \(f\cong f'\), if there are isomorphisms
\(\alpha:X\xrightarrow{\sim}X'\) and \(\beta:Y\xrightarrow{\sim}Y'\)
satisfying \(\beta\circ f=f'\circ\alpha\).
The following lemma shows that the multiplication map of \(V\),
\(Y_{Q,\boldsymbol\lambda}\) and \(H\) recovers
\(\rho_{Q,\boldsymbol\lambda}\) over \(I_{Q,\boldsymbol\lambda}^{\circ}\).

\begin{lemma}\label{lem:factorization-flags}
The multiplication map
\begin{equation}
  m:V\times Y_{Q,\boldsymbol\lambda}\times H\to U^-,\qquad(v,y,h)\mapsto vyh
\end{equation}
satisfies \(m\cong\id_{\A^{N-\dim R_Q}}\times
(\rho_{Q,\boldsymbol\lambda}|_{I_{Q,\boldsymbol\lambda}^{\circ}})\).
\end{lemma}

\begin{proof}
Recall that the map \(q:U^-\to R_Q\) is defined by \(q(g)_\gamma=g_{ij}g_{jj}^{-1}\) for each
arrow \(\gamma:j\to i\). For \(v\in V\) and \(g\in U^-\), we have \(q(vg)=v.q(g)\) because \(v\)
is block diagonal.
Since stability is \(G_W\)-equivariant,
\((q(g),vF_\bullet^{\boldsymbol\lambda})\in I_{Q,\boldsymbol\lambda}\)
if and only if
\((v^{-1}.q(g),F_\bullet^{\boldsymbol\lambda})\in I_{Q,\boldsymbol\lambda}\),
that is, if and only if \(q(v^{-1}g)=v^{-1}.q(g)\in Z\).
By \eqref{eq:quiver-chart}, we obtain an isomorphism
\begin{equation*}
U^-\times_{R_Q}I_{Q,\boldsymbol\lambda}^{\circ}
\xrightarrow{\sim}V\times q^{-1}(Z),\qquad
(g,(q(g),vF_\bullet^{\boldsymbol\lambda}))\longmapsto(v,v^{-1}g).
\end{equation*}
Moreover, the isomorphism \(\mu\) in (\ref{eq:mu}) restricting to \(Y_{Q,\boldsymbol\lambda}\times H\) gives an
isomorphism \(Y_{Q,\boldsymbol\lambda}\times H
\xrightarrow{\sim}q^{-1}(Z)\). Hence the map
\((v,y,h)\mapsto(vyh,\theta(v,y,h))\), where
\(\theta(v,y,h)=(q(vyh),vF_\bullet^{\boldsymbol\lambda})\),
identifies \(V\times Y_{Q,\boldsymbol\lambda}\times H\)
with the fiber product \(U^-\times_{R_Q}I_{Q,\boldsymbol\lambda}^{\circ}
\). Under this identification, we see that \(m\) is a projection to \(U^-\). Equivalently, the
commutative square
\begin{equation}\label{eq:family-identification}
\begin{tikzcd}[column sep=large,row sep=large]
V\times Y_{Q,\boldsymbol\lambda}\times H
\arrow[r,"\theta"]\arrow[d,"m"']
&I_{Q,\boldsymbol\lambda}^{\circ}
\arrow[d,"\rho_{Q,\boldsymbol\lambda}"]\\
U^-\arrow[r,"q"']&R_Q
\end{tikzcd}
\end{equation}
is cartesian.
Since \(q\) corresponds to the projection \(\A^{N-\dim R_Q}\times R_Q\to R_Q\)
under \eqref{eq:quiver-target-product}, the fiber product is
\(\A^{N-\dim R_Q}\times I_{Q,\boldsymbol\lambda}^{\circ}\), and \(m\) becomes
\(\id_{\A^{N-\dim R_Q}}\times(\rho_{Q,\boldsymbol\lambda}|_{I_{Q,\boldsymbol\lambda}^{\circ}})\).
\end{proof}

It remains to realize \(V\), \(Y_{Q,\boldsymbol\lambda}\) and \(H\) by
root-subgroup multiplication.
Let \(\alpha_1,\ldots,\alpha_{\dim W-1}\) be the simple roots of \(G=\SL(W)\),
in the standard order.
For \(1\le i\le j<\dim W\), write \([i,j]=\alpha_i+\cdots+\alpha_j\);
in particular, \([i,i]=\alpha_i\).
Any positive root of \(G\) has the form \([i,j]\).
For \(1\le j<i\le\dim W\), let \(E_{ij}\) be the matrix unit with \(1\) in the
position \((i,j)\) and \(0\) elsewhere.
Then \(U_{-[j,i-1]}=\{I_{\dim W}+tE_{ij}:t\in\C\}\).
As in the introduction, a sequence \(S=(\beta_1,\ldots,\beta_\ell)\) of
positive roots, now of arbitrary length \(\ell\), defines a
multiplication map \(\pi_S:\A^\ell\to U^-\).
For a set \(\Gamma\subseteq\Phi^+\) of positive roots of \(G\), let
\(S_\Gamma\) be a sequence obtained by ordering the roots
\([i,j]\)'s of \(\Gamma\) by nondecreasing \(i\).
Then we obtain the following result.
\begin{lemma}\label{lem:column-root-chart}
The map \(\pi_{S_\Gamma}\) is an
isomorphism onto \(I_{\dim W}+\bigoplus_{[j,i-1]\in\Gamma}\C E_{ij}\).
\end{lemma}

\begin{proof}
Write \(S_\Gamma=([j_1,i_1-1],\ldots,[j_\ell,i_\ell-1])\), where
\(\ell=|\Gamma|\) and \(j_1\le\cdots\le j_\ell\), and write
\(x_{-[j_k,i_k-1]}(t)=I_{\dim W}+c_ktE_{i_kj_k}\) with \(c_k\ne0\).
Then we get
\begin{equation*}
\begin{aligned}
\pi_{S_\Gamma}(t_1,\ldots,t_\ell)
&=x_{-[j_1,i_1-1]}(t_1)x_{-[j_2,i_2-1]}(t_2)\cdots x_{-[j_{\ell}, i_{\ell}-1]}(t_{\ell})\\
&=I_{\dim W}+\sum_{k=1}^{\ell}c_kt_kE_{i_kj_k}\\
&\quad+\sum_{m\ge2}\ \sum_{k_1<\cdots<k_m}
c_{k_1}\cdots c_{k_m}t_{k_1}\cdots t_{k_m}\,
E_{i_{k_1}j_{k_1}}\cdots E_{i_{k_m}j_{k_m}}.
\end{aligned}
\end{equation*}
For \(k<k'\), we have \(j_k\le j_{k'}<i_{k'}\), so \(j_k\ne i_{k'}\) and
\(E_{i_kj_k}E_{i_{k'}j_{k'}}=0\).
Hence the last sum vanishes, and
\begin{equation}\label{eq:linear-root-chart}
\pi_{S_\Gamma}(t_1,\ldots,t_\ell)
=I_{\dim W}+\sum_{k=1}^{\ell}c_kt_kE_{i_kj_k},
\end{equation}
which is an isomorphism.
\end{proof}

Recall that \(Q\) is an acyclic simple quiver. Let \(\boldsymbol\lambda\) be a flag type.
Applying the above lemma to \(V\),
\(Y_{Q,\boldsymbol\lambda}\) and \(H\), we obtain the following result.

\begin{lemma}\label{prop:quiver-realization}
There exists a sequence \(S_{Q,\boldsymbol\lambda}\) of positive
roots of \(G=\SL(W)\), with length
\(L=N+\sum_{a<b}\langle\lambda^a,\lambda^b\rangle_Q\), such that
\[
\pi_{S_{Q,\boldsymbol\lambda}}\cong\id_{\A^{N-\dim R_Q}}\times
(\rho_{Q,\boldsymbol\lambda}|_{I_{Q,\boldsymbol\lambda}^{\circ}}).
\]
\end{lemma}

\begin{proof}
Recall that \(H=G_W\cap U^-\), that \(V=\prod_{i\in Q_0}V_i\) with \(V_i\)
consisting of the block lower unitriangular matrices with diagonal blocks of
sizes \(\lambda_i^1,\ldots,\lambda_i^s\), and that
\(\mu^{-1}(Y_{Q,\boldsymbol\lambda})=Z\times R_Q'\times\{I_{\dim W}\}\), where \(Z\) is spanned by matrix units.
Hence each \(X\in\{V,Y_{Q,\boldsymbol\lambda},H\}\) is of the form
\(I_{\dim W}+\bigoplus_{[j,i-1]\in\Gamma_X}\C E_{ij}\) for some
\(\Gamma_X\subseteq\Phi^+\). Let \(\Gamma_{Y}=\Gamma_{Y_{Q,\boldsymbol{\lambda}}}\).
Set \(S_{Q,\boldsymbol\lambda}=S_{\Gamma_V}S_{\Gamma_{Y}}S_{\Gamma_H}\), the
concatenation of the three sequences.
By Lemma~\ref{lem:column-root-chart},
\(\pi_{S_{\Gamma_V}}\times\pi_{S_{\Gamma_Y}}\times\pi_{S_{\Gamma_H}}\) is an
isomorphism onto \(V\times Y_{Q,\boldsymbol\lambda}\times H\), and composing it
with \(m\) gives \(\pi_{S_{Q,\boldsymbol\lambda}}\).
Hence Lemma~\ref{lem:factorization-flags} gives the asserted isomorphism.
By \eqref{eq:quiver-dimension},
\(L=\dim V+\dim Y_{Q,\boldsymbol\lambda}+\dim H
=N-\dim R_Q+\dim I_{Q,\boldsymbol\lambda}
=N+\sum_{a<b}\langle\lambda^a,\lambda^b\rangle_Q\).
\end{proof}

We now apply the above lemma to the unicyclic incidence
projections of Lemma~\ref{lem:unicyclic-flag-projection}.

\begin{theorem}\label{cor:unicyclic-roots}
Let \(Q=(Q_0,Q_1)\) be a connected acyclic simple quiver such that \(|Q|\)
is unicyclic, and let \(n=\#(Q_0)\).
Let \(\boldsymbol{\vec{r}}=(r_1,\ldots,r_s)\) be a composition of \(r\). When
\(\dim W_i=r\) and \(\lambda_i^a=r_a\) for all \(i\in Q_0\) and \(1\le a\le s\),
the sequence \(S_{Q,\boldsymbol{\vec{r}}}:=S_{Q,\boldsymbol\lambda}\) in
Lemma~\ref{prop:quiver-realization} is a dominant root sequence of \(\SL(W)\), with length
\(\binom{nr}{2}\), and
\begin{equation}\label{eq:unicyclic-normal-form}
\pi_{S_{Q,\boldsymbol{\vec{r}}}}\sim_{\mathrm{bir}}
\id_{\A^{N-r^2}}\times p_{\boldsymbol{\vec{r}}}.
\end{equation}
Consequently,
\(\deg\pi_{S_{Q,\boldsymbol{\vec{r}}}}=r!/(r_1!\cdots r_s!)\), and for
\(\boldsymbol{\vec{r}}=(1,\ldots,1)\), the extension
\(\C(\A^N)/\pi_{S_{Q,\boldsymbol{\vec{r}}}}^*\C(U^-)\)
is Galois with
symmetric group \(\mathfrak S_r\).
\end{theorem}

\begin{proof}
Since \(\# (Q_0)=\#(Q_1)=n\), Lemma \ref{lem:quiver-incidence-dimension}
gives \(\dim I_{Q,\boldsymbol\lambda}
=\dim R_Q=nr^2\).
By Lemma~\ref{lem:unicyclic-flag-projection},
\(\rho_{Q,\boldsymbol\lambda}\) is dominant and
\(\rho_{Q,\boldsymbol\lambda}\sim_{\mathrm{bir}}
\id_{(\GL_r)^{n-1}}\times p_{\boldsymbol{\vec{r}}}\).
Since \(I_{Q,\boldsymbol\lambda}^{\circ}\) is dense in
\(I_{Q,\boldsymbol\lambda}\), Lemma~\ref{prop:quiver-realization}
implies that the length \(L=N\) and that \(\pi_{S_{Q,\boldsymbol{\vec{r}}}}\) is dominant.
By the fact that \((\GL_r)^{n-1}\) is rational, we obtain
\begin{equation*}
\pi_{S_{Q,\boldsymbol{\vec{r}}}}
\sim_{\mathrm{bir}}
\id_{\A^{N-nr^2}}\times\rho_{Q,\boldsymbol\lambda}
\sim_{\mathrm{bir}}
\id_{\A^{N-r^2}}\times p_{\boldsymbol{\vec{r}}}.
\end{equation*}
Generic degrees and Galois groups are invariant under \(\sim_{\mathrm{bir}}\)
and under products with affine spaces, so the remaining assertions follow from
Proposition~\ref{prop:flag-cover}.
\end{proof}

Therefore, Theorem~\ref{thm:flags} follows by applying
Theorem~\ref{cor:unicyclic-roots} to the triangle quiver
\begin{equation}\label{eq:tr-quiver}
Q:1\to2\to3,\ 1\to3.
\end{equation}
In fact, the sequence \(S_{Q,\boldsymbol{\vec{r}}}\) in
Theorem~\ref{cor:unicyclic-roots} is explicit: it is the concatenation
\(S_{\Gamma_V}S_{\Gamma_Y}S_{\Gamma_H}\), where \(\Gamma_V\), \(\Gamma_Y\) and
\(\Gamma_H\) are read off from the matrix forms of \(V\),
\(Y_{Q,\boldsymbol\lambda}\) and \(H\).
We illustrate an example of \(S_{Q,\boldsymbol{\vec{r}}}\) for the triangle quiver \(Q\) above and
\(\boldsymbol{\vec{r}}=(1,1)\).

\begin{example}\label{ex:A5}
Let \(U_2^-\subset\GL_2\) be the lower unitriangular
subgroup and \(\mathfrak b_2\subset\mathfrak{gl}_2\)
the upper triangular subalgebra. Then
\begin{equation*}
\begin{aligned}
V=H
&=\left\{
\begin{pmatrix}
u_1&0&0\\
0&u_2&0\\
0&0&u_3
\end{pmatrix}:u_i\in U_2^-
\right\},\qquad
Y_{Q,\boldsymbol\lambda}
&=\left\{
\begin{pmatrix}
I_2&0&0\\
B_{21}&I_2&0\\
B_{31}&B_{32}&I_2
\end{pmatrix}:B_{ij}\in\mathfrak b_2
\right\}.
\end{aligned}
\end{equation*}
Thus \(\Gamma_Y=\{[2l-1,2k-2],\,[2l,2k-2],\,[2l,2k-1]:1\le l<k\le3\}\),
and we obtain a possible sequence 
\(S_{\Gamma_{Y}}=\bigl(
[1,2],[1,4],
\alpha_2,[2,3],[2,4],[2,5],
[3,4],\alpha_4,[4,5]
\bigr).\)
Since \(S_{\Gamma_V}=S_{\Gamma_H}=(\alpha_1,\alpha_3,\alpha_5)\),
the concatenation \(S_{\Gamma_V}S_{\Gamma_Y}S_{\Gamma_H}\) is the sequence 
\[S=\bigl(
\alpha_1,\alpha_3,\alpha_5,
[1,2],[1,4],
\alpha_2,[2,3],[2,4],[2,5],
[3,4],\alpha_4,[4,5],
\alpha_1,\alpha_3,\alpha_5
\bigr).\]
By \eqref{eq:unicyclic-normal-form},
\(\pi_S\sim_{\mathrm{bir}}
\id_{\A^{11}}\times p_{(1,1)}\);
consequently, \(\deg\pi_S=2\).
\end{example}

\begin{remark}\label{rem:other-quivers}
If \(Q\) is acyclic and \(\rho_{Q,\boldsymbol\lambda}\) is birational, then
\(S_{Q,\boldsymbol\lambda}\) is a birational sequence in the sense of
\cite{FFL2017}.
Thus new birational sequences may be constructed from quiver
representations, by finding acyclic quivers \(Q\) and flag types
\(\boldsymbol\lambda\) for which \(\rho_{Q,\boldsymbol\lambda}\) is
birational.

Similarly, dominant incidence projections with
\(\dim I_{Q,\boldsymbol\lambda}=\dim R_Q\) yield dominant root
sequences of the same generic degree.
The triangle quiver in \eqref{eq:tr-quiver} is of affine type
\(\widetilde A_2\); it would be interesting to treat other affine quivers
in the same way.
\end{remark}

\section{Extension via Levi subgroups}\label{sec:levi}

In this section, we extend the type \(A\) examples to other types
via Levi subgroups.
The extensions preserve dominance and generic degree.

Let \(G\) be a connected semisimple algebraic group, with \(T\),
\(U^-\), \(\Phi\) and \(\Phi^+\)  as defined in the introduction.
Set \(N=|\Phi^+|=\dim U^-\), and write \(\Delta\subset\Phi^+\) for the set of
simple roots.
For \(I\subseteq\Delta\), let \(L_I\) be the corresponding
standard Levi subgroup, and set
\(\Psi=\Phi\cap\mathbb Z I\),
\(\Psi^+=\Psi\cap\Phi^+\),
\(U^-_\Psi=U^-\cap L_I\), and \(k=|\Psi^+|\).
The following lemma shows that extending a sequence in \(\Psi^+\) by the
roots of \(\Phi^+\setminus\Psi^+\) preserves dominance and generic degree.

\begin{lemma}\label{prop:levi}
Let \(S_\Psi=(\beta_1,\ldots,\beta_k)\) be a sequence
in \(\Psi^+\), with root-subgroup multiplication map
\(\pi_{S_\Psi}:\A^k\to U^-_\Psi\). Let
\(\gamma_1,\ldots,\gamma_{N-k}\) be any ordering of
\(\Phi^+\setminus\Psi^+\). Set
\(S=(\beta_1,\ldots,\beta_k,\gamma_1,\ldots,\gamma_{N-k})\).
Then
\begin{equation}\label{eq:levi-extension}
\pi_S\cong
\pi_{S_\Psi}\times\id_{\A^{N-k}}.
\end{equation}
Consequently, \(\pi_S\) is dominant if and only if 
\(\pi_{S_\Psi}\) is dominant. Moreover, when \(\deg\pi_S\) is finite, we have \(\deg\pi_S=\deg\pi_{S_\Psi}\).
\end{lemma}
\begin{proof}
By \cite[Theorem~21.68(a) and Remark~21.69]{Milne2017}, for any ordering of
the positive roots of a connected reductive group, multiplication of the
negative root subgroups is an isomorphism onto the negative maximal
unipotent subgroup.
Applying this to \(L_I\) and to \(G\), with the roots of \(\Psi^+\) ordered
first, we obtain an isomorphism
\begin{equation*}
\nu:U^-_\Psi\times\A^{N-k}\xrightarrow{\sim}U^-,\qquad
(u,\mathbf z)\longmapsto
u\,x_{-\gamma_1}(z_1)\cdots x_{-\gamma_{N-k}}(z_{N-k}).
\end{equation*}
Since \(\pi_S(\mathbf t,\mathbf z)=\nu(\pi_{S_\Psi}(\mathbf t),\mathbf z)\),
we get \(\pi_S=\nu\circ(\pi_{S_\Psi}\times\id_{\A^{N-k}})\), which gives
\eqref{eq:levi-extension}.
Taking the product with \(\id_{\A^{N-k}}\) preserves dominance and generic
degree, which gives the remaining assertions.
\end{proof}

\begin{corollary}\label{cor:classical}
For every integer \(d\ge2\), there exists a dominant
root sequence of generic degree \(d\) in type \(A_n\) for every \(n\ge3d-1\),
and in types \(B_n\), \(C_n\) and \(D_n\) for every \(n\ge3d\).
Moreover, there exists a dominant root sequence of generic degree \(2\) in
each of types \(E_6\), \(E_7\) and \(E_8\).
\end{corollary}

\begin{proof}
Taking \(r=d\) and \(\boldsymbol{\vec{r}}=(1,d-1)\) in
Theorem~\ref{thm:flags} gives a dominant sequence of degree \(d\) in type
\(A_{3d-1}\).
Recall that we have the inclusions of standard Levi root subsystems
\(A_{3d-1}\subseteq A_n\) for \(n\ge3d-1\),
\(A_{3d-1}\subseteq A_{n-1}\subseteq B_n,C_n,D_n\) for \(n\ge3d\),
and \(A_5\subseteq E_6,E_7,E_8\).
Hence, if \(G\) is of type \(A_n\) with \(n\ge3d-1\), of type \(B_n\), \(C_n\)
or \(D_n\) with \(n\ge3d\), or of type \(E_6\), \(E_7\) or \(E_8\) with \(d=2\),
then there exists \(I\subseteq\Delta\) such that \(\Psi\) is of type \(A_{3d-1}\).
The central isogeny \(\SL_{3d}\to[L_I,L_I]\) identifies the corresponding
root subgroups up to rescaling their parameters
\cite[Proposition~23.5]{Milne2017}, so the above sequence is a dominant
sequence in \(\Psi^+\) of degree \(d\).
Lemma~\ref{prop:levi} then extends it to a dominant sequence for \(G\) of the
same degree.
\end{proof}

\begin{remark}\label{rem:ranks}
The lower bounds \(n\ge3d-1\) and \(n\ge3d\) on the rank in
Corollary~\ref{cor:classical} are not sharp in general, since
Theorem~\ref{thm:flags} realizes every multinomial coefficient
\(r!/(r_1!\cdots r_s!)\) as a generic degree in type \(A_{3r-1}\).
For example, \(\boldsymbol{\vec{r}}=(1,1,1)\) gives degree \(6\) already in
type \(A_8\), rather than \(A_{17}\).
\end{remark}

\section*{Acknowledgements}
OpenAI's GPT-6 Pro was used in an author-directed 
process of mathematical exploration.
The authors independently generalized the constructions 
and developed the arguments. 
GPT-6 Astra and Anthropic's Claude (Opus 5.5) subsequently assisted 
with polishing the manuscript. The authors verified all mathematical 
statements and proofs and take full responsibility for the paper.

\end{document}